\documentclass[a4paper,12pt]{article}

\makeatletter
\@addtoreset{footnote}{page}
\makeatother

\renewenvironment{enumerate}
{\begin{list}{\parbox{2em}{(\arabic{enumi})}}{
\usecounter{enumi}
\setlength{\topsep}{0ex}
\setlength{\itemindent}{1em}
\setlength{\leftmargin}{2em} 
\setlength{\rightmargin}{1em}
\setlength{\labelsep}{1em}
\setlength{\labelwidth}{3em}%
\setlength{\itemsep}{0.1em}
\setlength{\parsep}{0em}
\setlength{\listparindent}{0em}
}}{\end{list}}

\usepackage{amsthm}
\usepackage{amsmath,amssymb,latexsym,amsfonts,mathrsfs}

\renewcommand{\d}{{\rm d}} 

\newcommand{\inde}{\perp \!\!\! \perp} 

\newcommand{\dist}{\stackrel{{\rm d}}{=}}

\newcommand{\tend}[2]{\mathrel{\mathop{\longrightarrow}\limits^{#1}_{#2}}}

\newcommand{\supre}[2]{\mathrel{\mathop{\sup}\limits^{#1}_{#2}}}

\newcommand{\maxi}[2]{\mathrel{\mathop{\max}\limits^{#1}_{#2}}}

\renewcommand{\hat}{\widehat}
\renewcommand{\tilde}{\widetilde}
\renewcommand{\bar}{\overline}

\newcommand{\absol}[1]{\left| #1 \right|} 
\newcommand{\norm}[1]{\left\| #1 \right\|} 
\newcommand{\rbra}[1]{\!\left( #1 \right)} 
\newcommand{\cbra}[1]{\!\left\{ #1 \right\}} 
\newcommand{\sbra}[1]{\!\left[ #1 \right]} 

\newcommand{\bE}{\ensuremath{\mathbb{E}}}

\newcommand{\bP}{\ensuremath{\mathbb{P}}}

\newcommand{\cB}{\ensuremath{\mathcal{B}}}

\newcommand{\cF}{\ensuremath{\mathcal{F}}}

\newcommand{\cT}{\ensuremath{\mathcal{T}}}

\newcommand{\cX}{\ensuremath{\mathcal{X}}}

\theoremstyle{plain}
\newtheorem{Thm}{Theorem}[section]

\newtheorem{Lem}[Thm]{Lemma}

\theoremstyle{definition}

\newtheorem{Def}[Thm]{Definition}

\newcommand{\Proof}[2][Proof]{\begin{proof}[{#1}] #2 \end{proof}}

\numberwithin{equation}{section}

\makeatletter
\renewcommand\section{\@startsection {section}{1}{\z@}%
                                   {-3.5ex \@plus -1ex \@minus -.2ex}%
                                   {2.3ex \@plus.2ex}%
                                   {\normalfont\large\bf}}
\makeatother

\makeatletter
\renewcommand\subsection{\@startsection {subsection}{1}{\z@}%
                                   {-3.5ex \@plus -1ex \@minus -.2ex}%
                                   {2.3ex \@plus.2ex}%
                                   {\normalfont\normalsize\bf}}
\makeatother
\begin{document}

\begin{center}
{\Large \bf 
Results for convergence rates associated with renewal processes
}
\end{center}
\footnotetext{
This research was supported by RIMS and by ISM.}
\begin{center}
Luis Iv\'an Hern\'andez Ru\'{\i}z \footnote{
Graduate School of Science, Kyoto University.}\footnote{
The research of this author was supported
by JSPS Open Partnership Joint Research Projects grant no. JPJSBP120209921.
}
\end{center}

\begin{abstract}
Via a coupling argument, it is proved that the solution to a renewal equation has a power law decay rate in the case of a spread out interarrival distribution. By the regenerative property, the convergence in distribution for the recurrence times of the renewal process is established as well as that of the compensator of the renewal counting process. All results are proved under mild assumptions of existence of moments.
\end{abstract}

\section{Introduction}
The Key Renewal Theorem states the limiting behaviour of the solution to the \emph{renewal equation}
\begin{align}
\label{RenewalEquation}
Z(t)=z(t)+\int_0^t Z(t-u)F(\d u),\quad t\ge0,
\end{align}
where $z:[0,\infty)\rightarrow[0,\infty)$ is known, $Z:[0,\infty)\rightarrow[0,\infty)$ is unknown and $F$ is a given probability distribution on $[0,\infty)$. If additionally $z$ is bounded on finite intervals, i.e. $\sup_{0\le t\le T}\absol{z(t)}<\infty$ for all $T<\infty$, then the only solution to \eqref{RenewalEquation} which is bounded on finite intervals is given by
\begin{align}
\label{RenewalEquationSolution}
Z(t)=\Phi*z(t)=\int_0^t z(t-u)\Phi(\d u),\quad t\ge0,
\end{align}
where the increasing function $\Phi(t)=\Phi\rbra{[0,t]}$ is given as 
\begin{align}
\label{RenewalFunction}
\Phi=\sum_{n\ge0}F^{*n},
\end{align}
with
\begin{align}
F^{*0}=\delta_0,\;\;\text{and }\; F^{*(n+1)}(t)=F^{*n}*F(t):=\int_0^tF^{*n}(t-s)F(\d s),\quad t\ge0.
\end{align}
The function $\Phi(t)$ is called the \emph{renewal function}, and its induced measure $\Phi(\d t)$ is called the \emph{renewal measure}.

Lindvall \cite{LindvallCoupling} found some of the first results on the rate of convergence between two differently delayed renewal measures with the aid of coupling techniques introduced by Athreya--Ney \cite{AthreyaNey}. More recently, Willmot--Cai--Lin \cite{WillmotLundberg} have found certain sharp estimates for the asymptotics of the solution \eqref{RenewalEquationSolution} under a variety of assumptions on $z$ and $F$. Asmussen--Foss--Korshunov \cite{asmussen2013asymptotics} considered this problem in the case of subexponential distributions. Yin--Zhao \cite{ChuancunNonExpAsym} studied non-exponential asymptotics in the case of \emph{defective distributions}, i.e. $F([0,\infty))<1$. Sgibnev \cite{Sgibnev1981} treated some cases where $Z(t)\rightarrow\infty$ as $t\rightarrow\infty$.

A different approach to the study of asymptotical properties in renewal theory concerns the limiting behaviour of processes constructed from a renewal process. Gakis--Sivazlian \cite{GakisRecurrence} obtained the limit distributions for the forward and backward recurrence times of a renewal process. Meanwhile, Godr\`eche--Luck \cite{godreche2000statistics} found results for occupation times of a renewal process using Laplace transform methods.

The purpose of this paper is twofold. First, we want to establish power law decay rates for the Key Renewal Theorem under the assumption of existence of moments. In particular, we are interested in the case of a \emph{spread out} interarrival distribution. 
\begin{Def}[\textbf{\cite[Sec VII. p.186]{Asmussen2003}}]
\label{DefSpreadOut}
A distribution $F$ on $[0,\infty)$ is called \emph{spread out} if for some $n\ge1$, there exists a nonnegative measure $G$ such that $0\neq G\le F^{*n}$ and $G$ is absolutely continuous with respect to the Lebesgue measure. The measure $G$ is called an \emph{absolutely continuous component} of $F^{*n}$.
\end{Def}
We obtain our result by means of a coupling argument as was done in the proof of the following Theorem (see for example \cite[Theorem VII.2.10]{Asmussen2003}).

\begin{Thm}[\textbf{Lund--Meyn--Tweedie \cite{LundExponential}}]
\label{ExpDecay}
Assume that the distribution $F$ is spread out with finite mean $m^{-1}:=\int_0^\infty xF(\d x)<\infty$ and that for some $\eta>0$, $\int_0^\infty e^{\eta x}F(\d x)<\infty$. Take $0<\epsilon<\eta$. If the function $z$ in \eqref{RenewalEquation} is measurable with $z(x)=O(e^{-\delta x})$ as $x\rightarrow\infty$ for some $\delta>\epsilon$, then
\begin{align}
\Phi*z(x)=m\int_0^\infty z(y)\d y+O(e^{-\epsilon x})\quad\text{as }x\rightarrow\infty.
\end{align}
\end{Thm}
Secondly, we want to study the convergence in distribution of the recurrence times and compensator of a renewal process as elements of the c\`adl\`ag space $D([0,1])$. We do so by noting that these processes evolve in cycles between renewals and satisfy the so-called \emph{regenerative property}.

Consider $\tau_1,\tau_2,\dots$ a sequence of i.i.d. random variables on $(0,\infty)$ defined on a probability space $(\Omega,\cF,\bP)$ with \emph{interarrival distribution} function $F(x)=\bP(\tau_1\le x)$ for $x\ge 0$. Additionally, consider a random variable $\tau_0$ on $[0,\infty)$ defined on the same space, which is independent of $\tau_1,\tau_2,\dots$, with \emph{delay distribution} $F_0(x)=\bP(\tau_0\le x)$, $x\ge0$, not necessarily equal to $F$. We will work in the setting of the following definition.
\begin{Def}
Let $\{\hat{S}_n\}_{n\ge0}$ be the sequence of partial sums given as
\begin{align}
\hat{S}_n=\tau_0+\tau_1+\dots+\tau_n,\quad n\ge0,
\end{align}
to which we associate the counting measure
\begin{align}
\hat{N}(B)=\sum_{n\ge0}1_{\{\hat{S}_n\in B\}},\quad B\in\cB([0,\infty)).
\end{align}
In particular, we identify the \emph{renewals} $\{\hat{S}_n\}_{n\ge0}$ with the \emph{delayed renewal process}
\begin{align}
\label{DefinitionRenewalProcess}
\hat{N}(t):=\hat{N}([0,t])=\sum_{n\ge0}1_{\{\hat{S}_n\le t\}}.
\end{align}
If instead we consider the partial sums
\begin{align}
S_0=0,\quad and\quad S_n=\tau_1+\dots+\tau_n,\quad n\ge1
\end{align}
then 
\begin{align}
N(t):=N([0,t])=\sum_{n\ge0}1_{\{S_n\le t\}}
\end{align}
is called a \emph{zero-delayed} or \emph{pure renewal process}.
\end{Def}

Notice that since the interarrival distributions do not have an atom at $x=0$, then necessarily $S_{n+1}>S_n$ on $\cbra{S_n<\infty}$ for all $n\ge0$, which implies that the renewal process is a \emph{simple point process}. From the definition of $N$, the following processes can be derived. 
\begin{Def}
\label{RecurrenceTimes}
Given a renewal process $N$, the \emph{backward recurrence time} $\cbra{A_t}_{t\ge0}$ and \emph{forward recurrence time} $\cbra{B_t}_{t\ge0}$ are given for $t\ge0$ as
\begin{align}
A_t=t-S_{N(t)-1},\;\quad B_t=S_{N(t)}-t.
\end{align}
\end{Def}

It is notable that $\cbra{A_t}_{t\ge0}$ and $\cbra{B_t}_{t\ge0}$ are time-homogeneous strong Markov processes (c.f. \cite[Proposition V.1.5]{Asmussen2003}). Denote by $(\cF^N_t)_{t\ge0}$ the augmentation of the natural filtration generated by the renewal process, given for each $t\ge0$ by $\sigma(N(s);0\le s\le t)$. The process $N$ is increasing and predictable, hence, by the Doob--Meyer decomposition theorem, it has an a.s. finite \emph{compensator} $\Lambda$ such that the process
\begin{align}
M:=N-\Lambda,
\end{align}
is an $(\cF^N_t)$-martingale.

Throughout this paper, we will make the following assumptions.
\begin{enumerate}
\item[\textbf{(A0)}]The interarrival distribution $F$ is spread out.
\item[\textbf{(A1)}]$F$ has a density $f$: $F(x)=\int_0^x f(s)\d s$.
\item[\textbf{(B0)}]$\tau$ has finite mean: $m^{-1}:=\bE\sbra{\tau}=\int_0^\infty xF(\d x)<\infty$.

\end{enumerate} 
Note that \textbf{(A1)} implies \textbf{(A0)}. We proceed now to state our main results.

The following Theorem improves the error term in Theorem \ref{ExpDecay}.
\begin{Thm}
\label{TheoremKRT}
Suppose that $F$ satisfies \textbf{\emph{(A0)}} and $\bE\sbra{\tau_1^s}=\int_0^\infty x^{s}\;F(\d x)<\infty$ for some $s\ge2$. Let $z:[0,\infty)\rightarrow[0,\infty)$ be a measurable function that is integrable, bounded, and $z(x)=O(x^{-r})$ as $x\rightarrow\infty$ for some $r>1$. Then, for $q$ such that $0\le q\le s-1$,
\begin{align}
\label{KRTDecay}
\Phi*z(x)=m\int_0^\infty z(y)\d y + O\rbra{x^{\max\cbra{1-r,-q}}}\quad\text{as }x\rightarrow\infty.
\end{align}
Moreover, if $z(x)=o(x^{-r})$ as $x\rightarrow\infty$, then \eqref{KRTDecay} holds  with $o\rbra{x^{\max\cbra{1-r,-q}}}$ instead of $O\rbra{x^{\max\cbra{1-r,-q}}}$.
\end{Thm}
The next result establishes the speed of convergence in distribution for the compensator of the renewal process.
\begin{Thm}
\label{MaxCompensator}
Let $N$ be a zero-delayed renewal process satisfying \textbf{(A1)} and \textbf{(B0)}. For any $p>0$, the convergence in distribution
\begin{align}
\rbra{\frac{1}{T^p}\sbra{\Lambda(Tv)-\Lambda(S_{N(Tv)-1})}}_{v\in[0,1]}\tend{d}{T\rightarrow\infty}0
\end{align}
holds in the Skorokhod topology.
\end{Thm}
Finally, we have the following Theorem that establishes the speed of convergence for the recurrence times.
\begin{Thm}
\label{MaxRecurrenceTime}
Let $N$ be a renewal process satisfying $\bE\sbra{\tau_1^s}=\int_0^\infty x^{s}\;F(\d x)<\infty$ for some $s\ge1$. Then, for any $p\le s$ the convergence in distribution
\begin{align}
\rbra{\frac{1}{T^{1/p}}\rbra{A_{Tv},B_{Tv}}}_{v\in[0,1]}\tend{d}{T\rightarrow\infty}(0,0)
\end{align}
holds in the Skorokhod topology.
\end{Thm}
The structure of the paper is the following. In Section \ref{SectionCoupling} we present in detail the construction of the coupling for renewal processes, and, in Section \ref{SectionDecay}, we use it to prove our Theorem \ref{TheoremKRT}. Section \ref{SectionRegenerative} consists of a brief explanation of regenerative processes. Section \ref{SectionConvergence} is dedicated to the proofs of our Theorems \ref{MaxCompensator} and \ref{MaxRecurrenceTime}. Finally, Section \ref{SecctionAppendixA} is an appendix with some results from the general theory.

\section{Coupling for renewal processes}
\label{SectionCoupling}
Throughout this section, we will suppose that \textbf{(A0)} and \textbf{(B0)} hold and we will review some known consequences of these assumptions. The first of such implications is that under \textbf{(B0)} we can make a delayed renewal process stationary, in the sense that for any $t>0$ the distribution of the increments $\cbra{\hat{N}(t+s)-\hat{N}(t)}_{s\ge0}$ does not depend on $t$, by choosing the delay distribution $\Pi=\pi(t)\d t$ with density $\pi$ given as
\begin{align}
\label{StationaryDistribution}
\pi:=m\bar{F}:=m(1-F),
\end{align}
in which case, the corresponding forward recurrence time $\hat{B}_t$ has common distribution $\Pi$ for all $t\ge0$ (see for example \cite[Proposition 4.2.I]{DVJ}).

Given a finite signed measure $\mu$ on the measurable space $(\Omega,\cF)$, denote by $\norm{\mu}_{\text{t.v.}}$ the total variation norm
\begin{align}
\norm{\mu}_{\text{t.v.}}:=\supre{}{B_1,B_2\in\cF}\rbra{\mu(B_1)-\mu(B_2)}.
\end{align}
Note that for measures, the total variation reduces to $\norm{\mu}_{\text{t.v.}}=\mu(\Omega)$.

Consider two stochastic processes $\cbra{X^\prime_t}_{t\ge0}$, $\cbra{X^{\prime\prime}_t}_{t\ge0}$, with the same state space and defined a priori on different probability spaces. By a \emph{coupling} of $X^{\prime}$, $X^{\prime\prime}$, we mean a pair $\rbra{\tilde{X}^{\prime},\tilde{X}^{\prime\prime}}$ and an associated random time $\cT$ (\emph{coupling time}), defined on a common probability space with
\begin{align}
\tilde{X}^{\prime}\dist X^{\prime} ,\quad\tilde{X}^{\prime\prime}\dist X^{\prime\prime},
\end{align}
and such that
\begin{align}
\tilde{X}^{\prime}_t=\tilde{X}^{\prime\prime}_t,\quad\text{for all }t\ge\cT.
\end{align}
We want to make use of the \emph{coupling inequality}. Let $\cbra{X^\prime_t}_{t\ge0}$, $\cbra{X^{\prime\prime}_t}_{t\ge0}$ be stochastic processes and $\cT\le\infty$ be a random time defined on a common probability space such that $X^\prime_t=X^{\prime\prime}_t$ for all $t\ge \cT$. Then,
\begin{align}
\label{CouplingInequality}
\norm{\bP\rbra{\theta_tX^\prime\in\cdot}-\bP\rbra{\theta_tX^{\prime\prime}\in\cdot}}_{\text{\emph{t.v.}}}\le 2\bP\rbra{\cT>t},
\end{align}
where $\theta_t$ stands for the shift operator given for any $t\ge0$ as
\begin{align}
\cbra{\theta_tX}_s=X_{t+s},
\end{align}
for any stochastic process $X$.

For two probability distributions $\lambda$ and  $\mu$, we define a measure $\lambda\wedge\mu$ as
\begin{align}
\lambda\wedge\mu(\cdot):=\int_\cdot(f\wedge g)\d(\lambda+\mu), 
\end{align} 
where $\lambda(\cdot)=\int_\cdot f\d(\lambda+\mu)$ and $\mu(\cdot)=\int_\cdot g\d(\lambda+\mu)$. Then we have
\begin{align}
\lambda\wedge\mu(\cdot)=\int_\cdot (\tilde{f}\wedge\tilde{g})\d\nu,
\end{align}
if $\lambda(\cdot)=\int_\cdot \tilde{f}\d\nu$ and $\mu(\cdot)=\int_\cdot \tilde{g}\d\nu$ for an arbitrary dominating measure $\nu$. Notice that
\begin{align}
\int(\tilde{f}-\tilde{g})\d \nu=1-1=0,
\end{align}
implies
\begin{align}
\int_{\cbra{\tilde{f}>\tilde{g}}}(\tilde{f}-\tilde{g})\d \nu=-\int_{\cbra{\tilde{f}\le\tilde{g}}}(\tilde{f}-\tilde{g})\d \nu.
\end{align}
 As a consequence,
\begin{align}
\norm{\lambda-\mu}_{\text{t.v.}}=&\int\absol{\tilde{f}-\tilde{g}}\d \nu
\label{}\\
=&\int_{\cbra{\tilde{f}>\tilde{g}}}(\tilde{f}-\tilde{g})\d \nu-\int_{\cbra{\tilde{f}\le\tilde{g}}}(\tilde{f}-\tilde{g})\d \nu
\label{}\\
=&2\int_{\cbra{\tilde{f}>\tilde{g}}}(\tilde{f}-\tilde{g})\d \nu
\label{}\\
=&2\int(\tilde{f}-\tilde{f}\wedge\tilde{g})\d \nu
\label{}\\
=&2\rbra{1-\norm{\lambda\wedge\mu}_{\text{t.v.}}}.
\label{TVIdentity1}
\end{align}
The following Lemma is sometimes referred to as \emph{maximal coupling}.
\begin{Lem}
\label{MaximalCoupling}
Given two probability distributions $F$, $G$ on a measurable space $(\cX,\cB(\cX))$, there exist random variables $X$, $Y$ defined on a common probability space $(\Omega,\cF,\bP)$ such that $X$ has distribution $F$, $Y$ has distribution $G$, and 
\begin{align}
\norm{F-G}_{\text{\emph{t.v.}}}=2\bP(X\neq Y).
\end{align}
\end{Lem}
\Proof{
Write $\delta:=\norm{F\wedge G}_{\text{t.v.}}$ and $H=\delta^{-1}(F\wedge G)$ and define the distributions
\begin{align}
F^\prime=\frac{F-\delta H}{1-\delta},\quad G^\prime=\frac{G-\delta H}{1-\delta}.
\end{align}
It is clear that
\begin{align}
F=\delta H+(1-\delta)F^\prime,\quad G=\delta H+(1-\delta)G^\prime.
\end{align}
Now we take independent random variables $\xi$, $X^\prime$, $Y^\prime$ and $Z$ defined on a common probability space $(\Omega,\cF,\bP)$ and such that $X^\prime$ has distribution $F^\prime$, $Y^\prime$ has distribution $G^\prime$, $Z$ has distribution $H$, and $\xi\sim\text{Bernoulli}(\delta)$. Define
\begin{align}
X:=\xi Z+(1-\xi)X^\prime,\;\;\text{and }\; Y:=\xi Z+(1-\xi)Y^\prime,
\end{align}
and notice that $X$ has distribution $F$, $Y$ has distribution $G$, and 
\begin{align}
\bP(X\neq Y)=\bP(\xi=0)=1-\delta=\frac{1}{2}\norm{F-G}_{\text{\emph{t.v.}}},
\end{align}
by \eqref{TVIdentity1}.
}
\subsection{Stone's decomposition}
\label{SubsectionStone}

It has been shown that when $F$ is spread out, the renewal measure can be written as a sum of a finite and an absolutely continuous component (see Stone \cite{Stone1966}). This decomposition is not unique, as it depends on the absolutely continuous component from Definition \eqref{DefSpreadOut}. For our purposes, it is convenient to select uniform components, which can always be done as shown in the following Lemma.
\begin{Lem}{\emph{\textbf{(see e.g., VII.1.2 from \cite{Asmussen2003})}}}
\label{RHPASpreadoutUniform}
If $F$ is spread out, then $F^{*n_0}$ has a uniform component on $(a,a+b)$ for some $a,b,n_0>0$.
\end{Lem}
\Proof{
Since $F$ is spread out, there exists a measure $G$ and an $n\ge0$ such that $0\neq G\le F^{*n}$ and $G$ has a density $g$ w.r.t. the Lebesgue measure.
Suppose that $g$ is bounded with compact support. Choose continuous functions $g_k\in L_1$  with compact supports such that
\begin{align}
\norm{g-g_k}_1=\int\absol{g(t)-g_k(t)}\d t\tend{}{k\rightarrow\infty}0.
\end{align}
Then
\begin{align}
\supre{}{\absol{x-x^\prime}<\delta}\absol{g_k*g(x)-g_k*g(x^\prime)}\le \supre{}{\absol{z-z^\prime}<\delta}\absol{g_k(z)-g_k(z^\prime)}\int g(y)\d y\tend{}{\delta\rightarrow0}0.
\end{align} 
So $g_k*g$ is uniformly continuous. Since
\begin{align}
\norm{g^{*2}-g*g_k}_\infty\le\norm{g}_\infty\norm{g-g_k}_1\tend{}{k\rightarrow\infty}0.
\end{align}
we see that $g^{*2}$ is continuous as a uniform limit of continuous functions, hence there exist $a,b,\delta>0$ such that $g^{*2}(x)\ge\delta$ for $x\in(a,a+b)$. Finally, take $n_0=2n$, and then
\begin{align}
G_0(\d x)=\delta 1_{(a,a+b)}(x)\d x
\end{align}
is a uniform component of $F^{*n_0}$.
}
We can then write Stone's decomposition for the renewal measure.
\begin{Thm}[\textbf{Stone \cite{Stone1966}}]
If the interarrival distribution $F$ is spread out, then we can write $\Phi=\Phi_1+\Phi_2$, where $\Phi_1$ and $\Phi_2$ are nonnegative measures on $[0,\infty)$, $\Phi_2$ is bounded (i.e. $\norm{\Phi_2}_{\text{\emph{t.v.}}}<\infty$) and $\Phi_1$ has a bounded density $\varphi_1(s)=\d\Phi_1(x)/\d x$ satisfying $\varphi_1(x)\tend{}{x\rightarrow\infty}m$.
\end{Thm}
The following proof is taken from Asmussen \cite{Asmussen2003} with a slight correction.
\Proof{
Let $G_0$ denote the uniform component of $F^{*n_0}$ which was given in the proof of Lemma \ref{RHPASpreadoutUniform}, and $g_0$ its density, given as
\begin{align}
g_0(x)=\frac{\norm{G_0}_{\text{t.v.}}}{b}1_{[a,a+b)}(x),\quad x\ge0,
\end{align}
and define $H:=F^{*n_0}-G_0$. We then note that
\begin{align}
\label{CyclicDecomp}
\Phi=\sum_{k=0}^{n_0-1}F^{*k}*\Phi_0,
\end{align}
where
\begin{align}
\Phi_0=&\sum_{n=0}^\infty F^{*nn_0}.
\end{align}
We can check by induction that
\begin{align}
F^{*nn_0}=&(H+G_0)^{*n}=G_0*\sum_{k=0}^{n-1}F^{*(n-k-1)n_0}*H^{*k}+H^{*n},
\end{align}
so $\Phi_0$ becomes
\begin{align}
\Phi_0=&\sum_{n=0}^{\infty}\rbra{G_0*\sum_{k=0}^{n-1}F^{*(n-k-1)n_0}*H^{*k}+H^{*n}}
\label{}\\
=&G_0*\sum_{k=0}^{\infty}H^{*k}*\sum_{n=k+1}^{\infty}F^{*(n-k-1)n_0}+\sum_{n=0}^{\infty}H^{*n}
\label{}\\
=&G_0*\sum_{n=0}^{\infty}H^{*n}*\Phi_0+\sum_{n=0}^{\infty}H^{*n}
\label{Phi0Decomp}
\end{align}
Then a decomposition of $\Phi_0$ is given as $\Phi_0=\Phi_0^{(1)}+\Phi_0^{(2)}$ where
\begin{align}
\Phi_0^{(2)}=\sum_{n=0}^{\infty}H^{*n},\quad \Phi_0^{(1)}=G_0*\Phi_0^{(2)}*\Phi_0.
\end{align}
Since $\norm{H}_{\text{t.v.}}=1-\norm{G_0}_{\text{t.v.}}<1$, we have that 
\begin{align}
\label{Normphi02}
\norm{\Phi_0^{(2)}}_{\text{t.v.}}=\frac{1}{1-\norm{H}_{\text{t.v.}}}=\frac{1}{\quad\norm{G_0}_{\text{t.v.}}}<\infty.
\end{align}
Moreover, since $G_0$ is absolutely continuous, so is $\Phi_0^{(1)}$, with density $\varphi_0^{(1)}=\Phi_0^{(2)}*(\Phi_0*g_0)$. From Blackwell's renewal theorem (c.f. \cite[Theorem 4.4.I]{DVJ}) we have that
\begin{align}
\label{Blackwellphi0g}
\Phi_0*g_0(x)=\frac{\norm{G_0}_\text{t.v.}}{b}\Phi_0((x-a-b,x-a])\tend{}{x\rightarrow\infty}\frac{m}{n_0}\norm{G_0}_\text{t.v.}.
\end{align}
From the subadditivity of the renewal function (c.f. \cite[Theorem V.2.4]{Asmussen2003}):
\begin{align}
 \Phi_0((x-a-b,x-a])\le\Phi_0(b)\quad\text{for all }x\ge0,
 \end{align}
and the total finiteness of $\Phi_0^{(2)}$, we have from the Dominated Convergence Theorem (DCT), \eqref{Normphi02} and \eqref{Blackwellphi0g} that
\begin{align}
\varphi_0^{(1)}(x)=\int_0^x\Phi_0*g_0(x-y)\Phi_0^{(2)}(\d y)\tend{}{x\rightarrow\infty}\frac{m}{n_0}\norm{G_0}_\text{t.v.}\norm{\Phi_0^{(2)}}_{\text{t.v.}}=\frac{m}{n_0}.
\end{align}
Finally, using \eqref{CyclicDecomp} and \eqref{Phi0Decomp}, we decompose $\Phi$ as $\Phi=\Phi_1+\Phi_2$ with
\begin{align}
\label{StonesDecomp}
\Phi_2=\sum_{k=0}^{n_0-1}F^{*k}*\Phi_0^{(2)},\quad\Phi_1=\sum_{k=0}^{n_0-1}F^{*k}*G_0*\Phi_0^{(2)}*\Phi_0=G_0*\Phi_0^{(2)}*\Phi,
\end{align} 
where
\begin{align}
\norm{\Phi_2}_{\text{t.v.}}=\frac{n_0}{\quad\norm{G_0}_{\text{t.v.}}}<\infty,
\end{align}
and $G_0*\Phi_0^{(2)}*\Phi$ has density $\varphi_1=\Phi_0^{(2)}*(\Phi*g_0)$ such that
\begin{align}
\varphi_1(x)=\int_0^x\Phi*g_0(x-y)\Phi_0^{(2)}(\d y)\tend{}{x\rightarrow\infty}m\norm{G_0}_\text{t.v.}\norm{\Phi_0^{(2)}}_{\text{t.v.}}=m.
\end{align}
Moreover, this density is also bounded since from subadditivity we have
\begin{align}
\Phi*g_0(x)=\frac{\norm{G_0}_\text{t.v.}}{b}\Phi([x-a-b,x-a))\le\frac{\norm{G_0}_\text{t.v.}}{b}\Phi(b),
\end{align}
from which we conclude that
\begin{align}
&\norm{\Phi*g_0}_{\infty}\le\frac{\norm{G_0}_\text{t.v.}}{b}\Phi(b)<\infty
\label{}\\
&\norm{\varphi_1}_{\infty}\le\norm{\Phi*g_0}_{\infty}\supre{}{x\ge0}\Phi_0^{(2)}((x-a-b,x-a])<\infty,
\end{align}
which concludes the proof.
}
In the remainder of this paper we use the decomposition given by \eqref{StonesDecomp}.
\subsection{Construction of the coupling}
We follow the construction of the coupling presented in \cite{thorisson2000coupling}.  First, we have the following Lemma for the forward recurrence time $B_t$ that has been introdued in Definition \ref{RecurrenceTimes}.
\begin{Lem}
\label{CommonUniformComponent}
For a zero-delayed spread out renewal process, there exist positive constants $b$ and $d$ such that the distributions of the $B_t$ with $t\ge d$ have a common uniform component on $(0,b)$. That is, for some $\delta\in(0,1)$ and all $t\ge d$,
\begin{align}
\bP(u<B_t\le b)\ge \delta\frac{v-u}{b},\quad 0<u<v<b.
\end{align}
\end{Lem}
We give a more detailed version of the proof found in \cite[Lemma VII.2.8]{Asmussen2003}.
\Proof{From Lemma \ref{RHPASpreadoutUniform} we have that for some $n_0\ge1$, $F^{*n_0}$ has a uniform component on an interval, so there exist constants $0<p<q$ and $\eta>0$ s.t. $F^{*n_0}(v)-F^{*n_0}(u)\ge \eta\rbra{v-u}$ for $0<p<u<v<q$. Let 
\begin{align}
b=\frac{(q-p)}{2},\quad\quad a=\frac{(p+q)}{2}.
\end{align}
When $p<z<a$, $0<u<v<b$, we have $\rbra{u+z,v+z}\subset\rbra{p,q}$, and hence from Lemma \ref{DistBt} we have
\begin{align}
\bP\rbra{u<B_t\le b}=&\int_0^t F\rbra{(t+u-s,t+b-s]}\Phi(\d s).
\label{UniformCompEq1}
\end{align}
Note that
\begin{align}
\Phi=\sum_{n\ge0}F^{*n}\ge\sum_{n\ge n_0-1}F^{*n}=F^{*(n_0-1)}*\Phi,
\label{}
\end{align}
from which we can say that
\begin{align}
&\bP\rbra{u<B_t\le b}
\label{}\\
=&\int_0^\infty 1_{[0,t]}(s_1) F\rbra{(t+u-s_1,t+b-s_1]}\Phi(\d s_1)
\label{}\\
=&\int_0^\infty\int_0^\infty 1_{(t+u,t+b]}(s_1+s_2)1_{[0,t]}(s_1) F(\d s_2)\Phi(\d s_1)
\label{}\\
\ge&\int_0^\infty\int_0^\infty 1_{(t+u,t+v]}(s_1+s_2)1_{[0,t]}(s_1) F(\d s_2)(F^{*(n_0-1)}*\Phi)(\d s_1)
\label{}\\
=&\int_0^\infty\int_0^\infty\int_0^\infty 1_{(t+u,t+v]}(s_2+s_3+s_4)1_{[0,t]}(s_3+s_4) F(\d s_2)F^{*(n_0-1)}(\d s_3)\Phi(\d s_4).
\label{}
\end{align}
It is easy to see that
\begin{align}
&1_{(t+u,t+v]}(s_2+s_3+s_4)1_{(t-p,t-a]}(s_4)
\label{}\\
\ge&1_{(t+u,t+v]}(s_2+s_3+s_4)1_{[0,t-a]}(s_4)1_{[0,v+p]}(s_3),
\label{}\\
\ge&1_{(t+u,t+v]}(s_2+s_3+s_4)1_{[0,t]}(s_3+s_4),
\end{align}
hence
\begin{align}
&\int_0^\infty\int_0^\infty\int_0^\infty 1_{(t+u,t+v]}(s_2+s_3+s_4)1_{[0,t]}(s_3+s_4) F(\d s_2)F^{*(n_0-1)}(\d s_3)\Phi(\d s_4)
\label{}\\
\ge&\int_0^\infty\int_0^\infty\int_0^\infty 1_{(t+u,t+v]}(s_2+s_3+s_4)1_{(t-p,t-a]}(s_4) F(\d s_2)F^{*(n_0-1)}(\d s_3)\Phi(\d s_4)
\label{}\\
=&\int_0^\infty\int_0^\infty 1_{(t+u,t+v]}(s_4+s_5)1_{(t-p,t-a]}(s_4) F^{*n_0}(\d s_5)\Phi(\d s_4)
\label{}\\
=&\int_{t-p}^{t-a}F^{*n_0}((t+u-s_4,t+v-s_4])\Phi(\d s_4).
\label{UniformCompEq2}
\end{align}
Using \eqref{UniformCompEq1} and \eqref{UniformCompEq2}, we obtain the inequality
\begin{align}
\bP\rbra{u<B_t\le b}\ge&\int_{t-p}^{t-a}F^{*n_0}((t+u-s,t+v-s])\Phi(\d s)
\label{}\\
\ge&\int_{p}^{a} F^{*n_0}\rbra{(u+z,v+z]}\Phi(t-\d z)
\label{}\\
\ge&\int_p^a\eta\rbra{v+z-u-z}\Phi(t-\d z)
\label{}\\
=&\eta(v-u)\sbra{\Phi(t-a)-\Phi(t-p)}.
\end{align}
Using Blackwell's renewal theorem, we have
\begin{align}
&\Phi(t-a)-\Phi(t-p)\tend{}{t\rightarrow\infty}m(a-p)=mb.
\end{align}
So, for $0<\delta<\min\cbra{1,m\eta b^2}$, there exists $d>0$ such that for all $t\ge d$,
\begin{align}
\Phi(t-a)-\Phi(t-p)\ge \frac{\delta}{b\eta},
\label{}
\end{align}
thus, for all $t\ge d$,
\begin{align}
\bP\rbra{u<B_t\le b}\ge \delta\frac{v-u}{b}.
\end{align}
The proof is complete.
}
We proceed with the construction of the coupling for the forward recurrence time of the renewal process.
\begin{Lem}[\textbf{Chapter 10, Theorem 6.2 from \cite{thorisson2000coupling}}]
\label{CouplingConstruction}
Let $\cbra{S_n}_{n\ge0}$ be a zero-delayed renewal process satisfying \textbf{(A0)} and \textbf{(B0)}, and $B=\cbra{B_t}_{t\ge0}$ its forward recurrence time. Let $\{\hat{S}_n\}_{n\ge0}$ be the stationary version of the renewal process whose forward recurrence time $\hat{B}=\{\hat{B}_t\}_{t\ge0}$ has stationary distribution $\Pi$. Then, the underlying probability space can be extended to support a coupling $(S^\prime,\hat{S^\prime})$ of $S$ and $\hat{S}$, and a geometric random variable $\sigma$ such that the coupling event occurs in a $\sigma$ number of trials. 
\end{Lem}
\Proof{
Let the positive constants $b$, $d$ and $\delta$ be as in Lemma \ref{CommonUniformComponent} so that   for all $t\ge d$, the distribution of $B_t$ has a common component $\delta U$ where $U$ has uniform distribution $\mu$ on $(0,b)$. Define a Markov process $(\eta_k,\hat{\eta}_k)_{k\ge0}$ in $[0,\infty)\times[0,\infty)$ in the following way:
\begin{align}
(\eta_0,\hat{\eta}_0):=(S_0,\hat{S}_0),\;\text{and } (\eta_k,\hat{\eta}_k):=(S_{N_{L_k-}},\hat{S}_{\hat{N}_{L_k-}}),\; k\ge1, 
\end{align}
where
\begin{align}
L_k:=\eta_k\vee\hat{\eta}_k+d.
\end{align}
Then, since $B$ is a time-homogeneous strong Markov process, we can see that conditionally on $(\eta_k,\hat{\eta}_k)=(s,\hat{s})$ the random variables
\begin{align}
\beta_{k+1}:=\eta_{k+1}-L_k,\quad\hat{\beta}_{k+1}:=\hat{\eta}_{k+1}-L_k,\quad k\ge0,
\end{align}
satisfy
\begin{align}
\bP(\beta_{k+1}\in A,\hat{\beta}_{k+1}\in \hat{A}\mid (\eta_k,\hat{\eta}_k)=(s,\hat{s}))=\bP(B_{((\hat{s}-s)_{+}+d)-}\in A)\bP(B_{((s-\hat{s})_{+}+d)-}\in \hat{A}).
\label{ForwardMarkov}
\end{align}
Thus, from Lemma \ref{CommonUniformComponent},
\begin{align}
\bP((\beta_{k+1},\hat{\beta}_{k+1})\in \cdot\mid (\eta_k,\hat{\eta}_k))\ge \delta^2 \mu\otimes \mu,\;k\ge0.
\end{align}
Using \cite[Chapter 3, Corollary 5.1]{thorisson2000coupling}, we can extend the underlying probability space to support 0-1 random variables $I_0,I_1,\dots$ such that for $k\ge0$,
\begin{align}
(S,\hat{S},I_0,I_1,\dots,I_{k-1})\inde I_k,\quad\text{given }{((\eta_k,\hat{\eta}_k),(\eta_{k+1},\hat{\eta}_{k+1}))},
\label{couplingcondind1}
\end{align}
and
\begin{align}
&\bP(I_k=1\mid\eta_k,\hat{\eta}_k)=\delta^2,
\label{Isindep}\\
&\bP((\beta_{k+1},\hat{\beta}_{k+1})\in \cdot\mid (\eta_k,\hat{\eta}_k)=(s,\hat{s}),I_k=1)=\mu\otimes \mu.
\end{align}
Fix an $m\ge1$. Because of \eqref{couplingcondind1} and the fact that $((\eta_i,\hat{\eta}_i),(\eta_{i+1},\hat{\eta}_{i+1}))$ is a measurable mapping of $(\eta_{k},\hat{\eta}_k)_{k=0}^{m}$, we have that for $0\le i<m$
\begin{align}
(S,\hat{S},I_0,I_1,\dots,I_{i-1})\inde I_i\quad\text{given }(\eta_{k},\hat{\eta}_k)_{k=0}^{m}.
\end{align}
In particular, for the $m$-tail of the Markov sequence it holds that
\begin{align}
(\eta_{m+k},\hat{\eta}_{m+k})_{k\ge 1}&\inde I_{i}\quad\text{given } \rbra{(\eta_{k},\hat{\eta}_k)_{k=0}^{m},I_0,\dots,I_{i-1}},
\end{align}
thus,
\begin{align}
(\eta_{m+k},\hat{\eta}_{m+k})_{k\ge 1}\inde& I_{m-1}\quad\text{given } \rbra{(\eta_{k},\hat{\eta}_k)_{k=0}^{m},I_0,\dots,I_{m-2}},
\label{}\\
(\eta_{m+k},\hat{\eta}_{m+k})_{k\ge 1}\inde& I_{m-2}\quad\text{given } \rbra{(\eta_{k},\hat{\eta}_k)_{k=0}^{m},I_0,\dots,I_{m-3}},
\label{}\\
\vdots&
\label{}\\
(\eta_{m+k},\hat{\eta}_{m+k})_{k\ge 1}\inde& I_0\quad\text{given } (\eta_{k},\hat{\eta}_k)_{k=0}^{m},
\end{align}
which combined with the Markov property of $(\eta_{k},\hat{\eta_{k}})_{k\ge 1}$ we can use to deduce that
\begin{align}
(\eta_{m+k},\hat{\eta}_{m+k})_{k\ge 1}\inde \rbra{(\eta_{k},\hat{\eta}_k)_{k=0}^{m},I_0,\dots,I_{m-1}}\quad\text{given }(\eta_m,\hat{\eta}_m).
\label{couplingcondind2}
\end{align}
Using \eqref{couplingcondind1} on $I_m$ we obtain that
\begin{align}
\rbra{(\eta_{k},\hat{\eta}_k)_{k=0}^{m},I_0,\dots,I_{m-1}}\inde I_m\quad \text{given }\rbra{(\eta_m,\hat{\eta}_m),(\eta_{m+k},\hat{\eta}_{m+k})_{k\ge 1}},
\end{align}
which in addition to \eqref{couplingcondind2} gives
\begin{align}
\rbra{(\eta_{k},\hat{\eta}_k)_{k=0}^{m},I_0,\dots,I_{m-1}}\inde \rbra{I_m,(\eta_{m+k},\hat{\eta}_{m+k})_{k\ge 1}}\quad \text{given }(\eta_m,\hat{\eta}_m).
\label{couplingcondind3}
\end{align}
Once again, by using \eqref{couplingcondind1} repeatedly we note that for each $n\ge1$
\begin{align}
(\eta_{k},\hat{\eta}_{k},I_k)_{k=0}^{m}\inde& I_{m+n}\quad\text{given } \rbra{(\eta_{m+k},\hat{\eta}_{m+k})_{k\ge1},I_{m+1},\dots,I_{m+n-1}},
\label{}\\
(\eta_{k},\hat{\eta}_{k},I_k)_{k=0}^{m}\inde& I_{m+n-1}\quad\text{given } \rbra{(\eta_{m+k},\hat{\eta}_{m+k})_{k\ge1},I_{m+1},\dots,I_{m+n-2}},
\label{}\\
\vdots&
\label{}\\
(\eta_{k},\hat{\eta}_{k},I_k)_{k=0}^{m}\inde& I_{m+1}\quad\text{given } (\eta_{m+k},\hat{\eta}_{m+k})_{k\ge1}.
\end{align}
We note from \eqref{couplingcondind3} that
\begin{align}
(\eta_{k},\hat{\eta}_{k},I_k)_{k=0}^{m}\inde(\eta_{m+k},\hat{\eta}_{m+k})_{k\ge 1}\quad \text{given }(\eta_m,\hat{\eta}_m,I_m).
\end{align}
We can then conclude that
\begin{align}
(\eta_{k},\hat{\eta}_{k},I_k)_{k=0}^{m}\inde (\eta_{m+k},\hat{\eta}_{m+k},I_{m+k})_{k\ge1}\quad\text{given } (\eta_m,\hat{\eta}_m,I_m),
\end{align}
hence, $(\eta_{k},\hat{\eta}_k,I_k)_{k\ge0}$ is a Markov process. Furthermore, from \eqref{Isindep} and \eqref{couplingcondind3} we know that the random variables $I_1,I_2,\dots$, are i.i.d. and $\bP(I_1=1)=\delta^2$. Thus,
\begin{align}
\sigma:=\inf\cbra{k\ge0: I_{k}=1},
\end{align}
is geometrically distributed with $\bP(\sigma=m)=(1-\delta^2)^m\delta^2$. Now, recall that $\eta_0=0$ (which implies $\eta_0\vee\hat{\eta}_0=\hat{\eta}_0$) and define
\begin{align}
\cT:=L_{\sigma}+U=T_\sigma,
\end{align}
where for $n\ge0$,
\begin{align}
T_n:=\hat{\eta_0}+d+\sum_{i=1}^{n}(\beta_i\vee\hat{\beta}_i+d)+U.
\label{DefinitionCouplingTime}
\end{align}
Let us define as well
\begin{align}
K:=N_{L_{\sigma}-},\;\hat{K}:=\hat{N}_{L_{\sigma}-}.
\end{align}
Since $\sigma$ is a stopping time with respect to the Markov process $(\eta_{k},\hat{\eta}_k,I_k)_{k\ge0}$, we obtain
\begin{align}
&\bP\rbra{(\beta_{\sigma+1},\hat{\beta}_{\sigma+1})\in\cdot\mid (\eta_{\sigma},\hat{\eta}_{\sigma})=(s,\hat{s})}
\label{}\\
=&\bP\rbra{(\beta_1,\hat{\beta}_1)\in\cdot\mid (\eta_{0},\hat{\eta}_{0})=(s,\hat{s}),I_0=1}=\mu\otimes \mu.
\label{SigmaIndep}
\end{align}
In other words, $\beta_{\sigma+1}$ and $\hat{\beta}_{\sigma+1}$ are i.i.d. with distribution $\mu$ and $(\beta_{\sigma+1},\hat{\beta}_{\sigma+1})$ is independent of $(\eta_{\sigma},\hat{\eta}_{\sigma})$ and thus of $L_{\sigma}$. From which we get that
\begin{align}
S_K=\eta_{\sigma+1}=L_{\sigma}+\beta_{\sigma+1}\dist L_{\sigma}+\hat{\beta}_{\sigma+1}=\hat{\eta}_{\sigma+1}=\hat{S}_K.
\end{align}
And due to $L_{\sigma}+\beta_{\sigma+1}\dist L_{\sigma}+U$, the processes $S$ and $\hat{S}$ have a common renewal at $\cT$. Taking $S^\prime=S$ and $\hat{S^\prime}$ with the same renewals as $\hat{S}$ before $\cT$ and with the same renewals as $S$ after $\cT$ gives the desired coupling.
}
\textbf{Remark.} Due to \eqref{Isindep} and the fact that $(\eta_k,\hat{\eta}_k,I_k)_{k\ge0}$ is a Markov process, we have for $i\in\{0,1\}$ and $k\ge1$ that
\begin{align}
&\bP\rbra{(\eta_{0},\hat{\eta}_{0})\in\cdot,I_0=i)}=\bP\rbra{I_0=i}\bP\rbra{(\eta_{0},\hat{\eta}_{0})\in\cdot},
\label{}
\end{align}
and
\begin{align}
&\bP\rbra{(\eta_{k},\hat{\eta}_{k})\in\cdot,I_k=i\mid (\eta_{k-1},\hat{\eta}_{k-1})=(s,\hat{s}),I_k=0}
\label{}\\
=&\bP\rbra{I_{k}=i}\bP\rbra{(\eta_{k},\hat{\eta}_{k})\in\cdot\mid (\eta_{k-1},\hat{\eta}_{k-1})=(s,\hat{s})}.
\end{align}
Since the event
\begin{align}
\cbra{\sigma=m}=\cbra{I_0=0,I_1=0,\dots,I_{m-1}=0,I_m=1},
\end{align}
we have that for any Borel sets $D_0, D_1,\dots, D_{m-1}$,
\begin{align}
&\bP\rbra{(\eta_{0},\hat{\eta}_{0})\in D_0,\dots,(\eta_{m-1},\hat{\eta}_{m-1})\in D_{m-1},\sigma=m}
\label{}\\
=&\int_{D_0}\cdots\int_{D_{m-1}}\bP\rbra{(\eta_{m-1},\hat{\eta}_{m-1})\in\d u_{m-1},I_{m-1}=0\mid (\eta_{m-2},\hat{\eta}_{m-2})=u_{m-2},I_{m-2}=0}\cdots
\label{}\\
\cdots&\bP\rbra{(\eta_{1},\hat{\eta}_{1})\in\d u_1,I_{1}=0\mid (\eta_{0},\hat{\eta}_{0})=u_0,I_{0}=0}\bP\rbra{(\eta_{0},\hat{\eta}_{0})\in\d u_0,I_0=0}
\label{}\\
=&\bP(I_m=1)\bP(I_{m-1}=0)\cdots\bP(I_0=0)\bP\rbra{(\eta_{0},\hat{\eta}_{0})\in D_0,\dots,(\eta_{m-1},\hat{\eta}_{m-1})\in D_{m-1}},
\end{align}
hence,
\begin{align}
&\bP\rbra{(\eta_{0},\hat{\eta}_{0})\in\cdot,\dots,(\eta_{m-1},\hat{\eta}_{m-1})\in\cdot,\sigma=m}
\label{}\\
=&\bP(\sigma=m)\bP\rbra{(\eta_{0},\hat{\eta}_{0})\in\cdot,\dots,(\eta_{m-1},\hat{\eta}_{m-1})\in\cdot}.
\label{SigmaIndepFromT}
\end{align}

\section{Decay rates for the Key Renewal Theorem}
\label{SectionDecay}
We work in the same context as the previous section. For the coupling time $\cT$ constructed in Lemma \ref{CouplingConstruction} we have the following result that we prove along the lines of \cite{LindvallCoupling}.

\begin{Lem}
\label{PowerCoupTime} 
Assume $\bE\sbra{\tau_1^s}<\infty$ for some $s\ge2$. Then $\bE\sbra{\cT^{q}}<\infty$ for $0< q\le s-1$.
\end{Lem}
\Proof{It suffices to show $\bE\sbra{\cT^{q}}<\infty$ for $q=s-1$. Let us find some estimates for $\bE\sbra{B_t^q}$. Consider first
\begin{align}
\bE\sbra{B_t}=&\bE\sbra{S_{N(t)}}-t
\label{}\\
=&m^{-1}\bE\sbra{N(t)}-t
\label{}\\
=&m^{-1}\Phi(t)-t.
\label{ExpctationBt}
\end{align}
Since $\Phi(t)/t\tend{}{t\rightarrow\infty}m$, for any $\epsilon_0>0$ we can find a $t_0>0$ such that 
\begin{align}
\Phi(t)<&(\epsilon_0+m)t,\quad\text{for all }t\ge t_0,
\label{}\\
\Phi(t)\le&\Phi(t_0)+(\epsilon_0+m)t,\quad\text{for all }t\ge 0.
\label{BoundExpectationBt}
\end{align}
Hence, for any arbitrary $\rho>0$, there exists a constant $a_\rho$ such that
\begin{align}
\bE\sbra{B_t}\le a_\rho+\rho t,\quad\text{for all }t\ge0.
\label{OvershootAlternativeBound}
\end{align} 
For $\bE\sbra{B_t^q}$, we can use Lemma \ref{DistBt} to find
\begin{align}
\bE\sbra{B_t^q}=&\int_0^t \int_0^\infty x^qF\rbra{t-u+\d x}\Phi(\d u)
\label{}\\
\le&\int_0^t \int_{t-u}^\infty x^qF\rbra{\d x}\Phi(\d u)
\label{}\\
\le&\rbra{\int_{0}^\infty x^qF\rbra{\d x}}\Phi(t),
\end{align}
where $\int_{0}^\infty x^qF\rbra{\d x}<\infty$. Then, using \eqref{BoundExpectationBt} with $\epsilon_0=1$, we find positive constants $a_1$ and $b_1$ independent of $k$ such that
\begin{align}
\bE\sbra{B_t^q}\le a_1+b_1 t,\quad t\ge0
\end{align}
and by Fatou's Lemma also
\begin{align}
\bE\sbra{B_{t-}^q}\le\liminf_{n\rightarrow\infty}\bE\sbra{B^q_{t-1/n}}\le a_1+b_1 t,\quad t\ge0.
\label{BoundOvershoot}
\end{align}

Since there exist constants $a_2$ and $b_2$ such that $(d+x)^q\le a_2+b_2x^q$, for all $x\ge0$, we have
\begin{align}
\bE\sbra{(d+\beta_k\vee\hat{\beta}_k)^q}\le a_0+b_0\bE\sbra{(\beta_k\vee\hat{\beta}_k)^q}.
\label{OvershootFirstEstimate}
\end{align}
From \eqref{ForwardMarkov} we have
\begin{align}
&\bE\sbra{\rbra{\beta_k\vee\hat{\beta}_k}^q\bigg| (\eta_{k-1},\hat{\eta}_{k-1})=(y,\hat{y})}
\label{}\\
\le&\bE\sbra{\beta_k^q+\hat{\beta}_k^q\bigg| (\eta_{k-1},\hat{\eta}_{k-1})=(y,\hat{y})}
\label{}\\
=&\bE\sbra{B_{((\hat{y}-y)_{+}+d)-}^q}+\bE\sbra{B_{((y-\hat{y})_{+}+d)-}^q}.
\label{TwoOvershoots}
\end{align}
Then, from \eqref{BoundOvershoot}, we obtain
\begin{align}
\bE\sbra{(\beta_k\vee\hat{\beta}_k)^q\mid \eta_{k-1},\hat{\eta}_{k-1}}
\le& 2a_1+2b_1 d+b_1\absol{\eta_{k-1}-\hat{\eta}_{k-1}}
\label{}\\
=& 2a_1+2b_1 d+b_1\absol{\beta_{k-1}-\hat{\beta}_{k-1}}
\label{}\\
\le& 2a_1+2b_1 d+b_1(\beta_{k-1}\vee\hat{\beta}_{k-1}).
\label{}
\end{align}
We define $a_3:=2(a_1+b_1d)$ and $b_3:=b_1$ to write
\begin{align}
\bE\sbra{(\beta_k\vee\hat{\beta}_k)^q\mid \eta_{k-1},\hat{\eta}_{k-1}}\le a_3+b_3(\beta_{k-1}\vee\hat{\beta}_{k-1}).
\label{OvershootRecurrence}
\end{align}
Now we shift our attention to $\bE\sbra{\beta_k\vee\hat{\beta}_k}$ for $k\ge1$. We follow the same steps as in \eqref{TwoOvershoots} (but with $q=1$), and use \eqref{OvershootAlternativeBound} with $\rho=1/2$ to find a constant $a_\rho$ such that
\begin{align}
\bE\sbra{\beta_k\vee\hat{\beta}_k\mid \eta_{k-1},\hat{\eta}_{k-1}}\le 2a_\rho+d+\frac{1}{2}\rbra{\beta_{k-1}\vee\hat{\beta}_{k-1}}.
\end{align}
By taking conditional expectation given $(\eta_{k-2},\hat{\eta}_{k-2}),(\eta_{k-3},\hat{\eta}_{k-3}),\dots$, recursively and using \eqref{OvershootAlternativeBound} with $\rho=1/2$ we ultimately obtain
\begin{align}
\bE\sbra{\beta_k\vee\hat{\beta}_k\mid \eta_{0},\hat{\eta}_{0}}\le a_\rho(2+1+\dots+2^{2-k})+d(1+2^{-1}+\dots+2^{1-k})+2^{-k}\rbra{\eta_{0}\vee\hat{\eta}_{0}},
\end{align}
which can be simplified by recalling that $\eta_0=0$ and by noting that for any $k\ge1$: $2+1+\dots+2^{2-k}< 4$, $1+2^{-1}+\dots+2^{1-k}< 2$, and $2^{-k}< 1$. In summary,
\begin{align}
\bE\sbra{\beta_k\vee\hat{\beta}_k\mid \eta_{0},\hat{\eta}_{0}}\le 4a_\rho+2d+\hat{\eta}_{0}=a_4+\hat{\eta}_{0},
\label{ExpectationBound}
\end{align}
where $a_4:=4a_\rho+2d$. Now note that since $\hat{\eta}_0$ is taken from the stationary distribution $\Pi$ and $0< q\le s-1$, then
\begin{align}
\bE\sbra{\hat{\eta}_0^{q}}=&\int_0^\infty x^{q}m\bar{F}(x)\d x
\label{}\\
=&\frac{m}{q+1}\int_0^\infty y^{q+1}F(\d y)<\infty.
\end{align}

There are two notable cases. Let us focus first on the case $q=s-1\ge1$. In this case $c_0:=\bE\sbra{(d+\hat{\eta_0})^q}<\infty$. We take the expectation on both sides of \eqref{ExpectationBound}  and \eqref{OvershootFirstEstimate}, which together with \eqref{OvershootRecurrence} yields
\begin{align}
\bE\sbra{(d+\beta_k\vee\hat{\beta}_k)^q}\le& a_2+b_2\rbra{a_3+b_3\bE\sbra{(\beta_k\vee\hat{\beta}_k)^q} }
\label{}\\
\le&a_2+b_2\cbra{a_3+b_3\rbra{a_4+c_0}}.
\label{OvershootCommonBound}
\end{align}
We take $a_q:=a_2+b_2\cbra{a_3+b_3\rbra{a_4+c_0}}$. Because $q\ge1$, we can apply Minkowski's inequality to obtain that
\begin{align}
\bE\sbra{T_n^{q}}^{1/q}\le& \bE\sbra{(d+\hat{\eta}_0)^q}^{1/q}+\sum_{k=1}^n\bE\sbra{\rbra{d+\beta_k\vee\hat{\beta}_k}^{q}}^{1/q}+(b^q)^{1/q}
\label{}\\
\le& c_0^{1/q}+na_q^{1/q}+b.
\end{align}
Since $\cbra{\sigma=m}$ is independent form $T_n$ by \eqref{SigmaIndepFromT} and \eqref{DefinitionCouplingTime}, and since $\sigma$ has geometric distribution, we have
\begin{align}
\bE\sbra{\cT^q}=&\bE\sbra{T_{\sigma}^q}
\label{}\\
=&\sum_{n=0}^\infty\bE\sbra{T_{n}^q}\bP\rbra{\sigma=n}
\label{}\\
\le&\sum_{n=0}^\infty(c_0^{1/q}+na_q^{1/q}+b)^q\delta^2(1-\delta^2)^n<\infty,
\end{align}
because $(1-\delta^2)^n$ vanishes exponentially.

For the case where $q=s-1<1$ we use the subadditivity of the function $x\mapsto x^q$, $x\ge0$, to get
\begin{align}
\bE\sbra{\cT^q}=&\bE\sbra{\rbra{\sum_{k=0}^n(d+\beta_k\vee\hat{\beta}_k)+U}^q;\sigma=n}
\label{}\\
\le&\bE\sbra{(d+\hat{\eta}_0)^q}+\sum_{k=1}^\infty\bE\sbra{ (d+\beta_k\vee\hat{\beta}_k)^q 1_{\cbra{\sigma\ge k}}}+b^q
\label{}\\
\le&c_0+\sum_{k=1}^\infty\bE\sbra{ \bE\sbra{(d+\beta_k\vee\hat{\beta}_k)^q 1_{\cbra{\sigma\ge k}}\mid \eta_0,\hat{\eta}_0}}+b^q.
\end{align}
Since $0<q<1$, we can apply Hölder's inequality with $p>1$ such that $q+\frac{1}{p}=1$ to obtain
\begin{align}
\bE\sbra{ \bE\sbra{(d+\beta_k\vee\hat{\beta}_k)^q 1_{\cbra{\sigma\ge k}}\mid \eta_0,\hat{\eta}_0}}\le&\bE\sbra{ \bE\sbra{(d+\beta_k\vee\hat{\beta}_k)\mid \eta_0,\hat{\eta}_0}^q\bE\sbra{ 1_{\cbra{\sigma\ge k}}\mid \eta_0,\hat{\eta}_0}^{1/p}}
\label{}\\
\le&\bE\sbra{ d^q+ \bE\sbra{\beta_k\vee\hat{\beta}_k\mid \eta_0,\hat{\eta}_0}^q}\bP\rbra{\sigma\ge k}^{1/p}.
\end{align}
From \eqref{ExpectationBound} we know that
\begin{align}
\bE\sbra{\beta_k\vee\hat{\beta}_k\mid \eta_0,\hat{\eta}_0}^q\le a_4^q+\hat{\eta}_0^q,
\end{align}
hence,
\begin{align}
\bE\sbra{ \bE\sbra{(d+\beta_k\vee\hat{\beta}_k)^q 1_{\cbra{\sigma\ge k}}\mid \eta_0,\hat{\eta}_0}}\le& \rbra{d^q+a_4^q+\bE\sbra{\hat{\eta}_0^q}}\bP\rbra{\sigma\ge k}^{1/p}
\label{}\\
=&c_q\bP\rbra{\sigma\ge k}^{1/p},
\end{align}
where $c_q:=d^q+a_4^q+\bE\sbra{\hat{\eta}_0^q}<\infty$. The random variable $\sigma$ is geometrically distributed so $\bP\rbra{\sigma\ge k}=(1-\delta^2)^k$. From the above, we have that
\begin{align}
\bE\sbra{\cT^q}\le& c_0+b^q+\sum_{k=1}^\infty c_q(1-\delta^2)^{k/p}<\infty,
\end{align}
since $(1-\delta^2)^{k/p}$ decreases exponentially in $k$. This ultimately means that $\bE\sbra{\cT^q}<\infty$ for $0<q\le s-1$, which concludes the proof.
}

Before proceeding, let us make the following observation. Let $h:[0,\infty)\rightarrow[0,\infty)$ be a non-decreasing function with $\lim_{t\rightarrow\infty}h(t)=\infty$ and suppose that $\bE\sbra{h(\cT)}<\infty$ for the coupling time $\cT$, then, by the DCT,
\begin{align}
h(t)\bP(\cT>t)\le&\bE\sbra{h(\cT)1_{\cbra{\cT>t}}}\tend{}{t\rightarrow\infty}0
\label{}
\end{align}
which means that
\begin{align}
\bP(\cT>t)=o\rbra{\frac{1}{h(t)}}\quad\text{as }t\rightarrow\infty.
\label{CouplingInequality2}
\end{align}
We use the previous remark to prove the following Lemma.
\begin{Lem}
\label{DecayRatesLemma}
Let $\{S_n\}_{n\ge0}$ be a zero-delayed renewal process that satisfies $\bE\sbra{\tau_1^s}<\infty$ for some $s\ge2$. Then, for $q$ such that $0< q\le s-1$,
\begin{align}
\norm{\bP\rbra{B_t\in\cdot}-\Pi(\cdot)}_{\text{\emph{t.v.}}}=o\rbra{t^{-q}}\; \text{as } t\rightarrow\infty,
\end{align}
where $\Pi$ is the stationary distribution given in \eqref{StationaryDistribution}, and
\begin{align}
\Phi_2([x,\infty))=o\rbra{x^{-q}}\; \text{as } x\rightarrow\infty, \quad\varphi_1(x)=m+o\rbra{x^{-q}} \;\text{as } x\rightarrow\infty.
\end{align}
\end{Lem}
\Proof{
Let $q$ such that $0<q<s-1$. Notice that the function $h:[0,\infty)\rightarrow[0,\infty)$ given by $h(t)=t^q$ is non-decreasing and $\bE\sbra{\cT^q}<\infty$, by Lemma \ref{PowerCoupTime}. By definition of $\Pi$ and $\cT$ we have
\begin{align}
\norm{\bP\rbra{B_t\in\cdot}-\Pi(\cdot)}_{\text{\emph{t.v.}}}=\norm{\bP\rbra{B_t\in\cdot}-\bP(\hat{B}_t\in\cdot)}_{\text{\emph{t.v.}}}\le2\bP\rbra{\cT>t},
\end{align}
and thus, from \eqref{CouplingInequality2},
\begin{align}
\lim_{t\rightarrow\infty}t^q\norm{\bP\rbra{B_t\in\cdot}-\Pi(\cdot)}_{\text{\emph{t.v.}}}\le\lim_{t\rightarrow\infty}2t^q\bP\rbra{\cT>t}=0,
\end{align}
which implies
\begin{align}
\label{DecayTVNorm}
&\norm{\bP\rbra{B_t\in\cdot}-\Pi(\cdot)}_{\text{\emph{t.v.}}}=o\rbra{t^{-q}}\quad\text{as }t\rightarrow\infty.
\end{align}
Now consider Stone's decomposition $\Phi=\Phi_1+\Phi_2$ obtained from the uniform component $G_0$ with density $g_0(x)=\frac{\norm{G_0}_{\text{t.v.}}}{b}1_{\cbra{a\le x<a+b}}$, $x\ge0$. Since $H=F^{*n_0}-G_0<F$, from the finiteness of $\bE\sbra{\tau_1^s}$, we get that 
\begin{align}
&\int_{0}^{\infty}x^{s}H(\d x)<\infty,
\end{align}
and hence,
\begin{align}
&\int_{0}^{\infty}x^s\Phi_2(\d x)<\infty,
\label{}
\end{align}
and
\begin{align}
\limsup_{x\rightarrow\infty}x^q\int_{[x,\infty)}\Phi_2(\d t)\le\lim_{x\rightarrow\infty}\int_{[x,\infty)}t^q\Phi_2(\d t)=0,
\end{align}
from which we get that
\begin{align}
\Phi_2([x,\infty))=o\rbra{x^{-q}}\; \text{as } x\rightarrow\infty.
\end{align}
From \cite[Sec. V. Theorem 2.4 (iii)]{Asmussen2003} we have that for a general renewal process
\begin{align}
&\Phi\rbra{(x-a-b,x-a]}
\label{}\\
=&\bE\sbra{N_{(x-a-b)+b}-N_{x-a-b}}
\label{}\\
=&\int_0^b\Phi(b-\xi)\bP\rbra{B_{x-a-b}\in\d\xi}
\label{}\\
=&\int_0^b\Phi(b-\xi)\Pi\rbra{\d\xi}+\int_0^b\Phi(b-\xi)\sbra{\bP\rbra{B_{x-a-b}\in\d\xi}-\Pi\rbra{\d\xi}}.
\label{}
\end{align}
Since
\begin{align}
\absol{\int_0^b\Phi(b-\xi)\sbra{\bP\rbra{B_{x-a-b}\in\d\xi}-\Pi\rbra{\d\xi}}}\le \Phi(b)\norm{\bP\rbra{B_t\in\cdot}-\Pi(\cdot)}_{\text{\emph{t.v.}}},
\end{align}
we can use \eqref{DecayTVNorm} and Lemma \ref{RHPALinearSolution} to obtain,
\begin{align}
&\int_0^b\Phi(b-\xi)\Pi\rbra{\d\xi}+\int_0^b\Phi(b-\xi)\sbra{\bP\rbra{B_{x-a-b}\in\d\xi}-\Pi\rbra{\d\xi}}
\label{}\\
=&mb+o\rbra{\rbra{x-a-b}^{-q}}
\label{}\\
=&mb+o\rbra{x^{-q}}.
\label{}
\end{align}
In other words, $\Phi\rbra{(x-a-b,x-a]}=mb+o\rbra{x^{-q}}$ as $x\rightarrow\infty$. From Stone's decomposition the density $\varphi_1$ is given as
\begin{align}
\varphi_1(x)=&\int_0^x\Phi*g_0(x-y)\Phi^{(2)}_0(\d y),
\end{align}
where from the previous reasoning
\begin{align}
\Phi*g_0(x)-m\norm{G_0}_{\text{t.v.}}=o(x^{-q})
\end{align}
Notice that
\begin{align}
&\absol{\int_{0}^{x/2}\rbra{\Phi*g_0(x-y)-m\norm{G_0}_{\text{t.v.}}}\Phi_0^{(2)}(\d y)}
\label{}\\
\le&\supre{}{y^\prime\ge x/2}\absol{\Phi*g_0(x-y^\prime)-m\norm{G_0}_{\text{t.v.}}}\norm{\Phi_0^{(2)}}_{\text{t.v.}}=o(x^{-q}),
\end{align}
and
\begin{align}
&\absol{\int_{x/2}^{x}\rbra{\Phi*g_0(x-y)-m\norm{G_0}_{\text{t.v.}}}\Phi_0^{(2)}(\d y)}
\label{}\\
\le&\rbra{\norm{\Phi*g_0}_{\infty}+m\norm{G_0}_{\text{t.v.}}}\Phi_0^{(2)}((x/2,\infty))=o(x^{-q}),
\end{align}
from which we conclude that
\begin{align}
\varphi_1(x)=&m\norm{G_0}_{\text{t.v.}}\norm{\Phi_0^{(2)}}_{\text{t.v.}}+o(x^{-q})
\label{}\\
=&m+o(x^{-q}),
\end{align}
by \eqref{Normphi02}.
}
We proceed to the proof of Theorem \ref{TheoremKRT}.

\Proof[Proof of Theorem \ref{TheoremKRT}]{We prove the case where $z(x)=o(x^{-r})$ as $x\rightarrow\infty$ and the case for $z(x)=O(x^{-r})$ as $x\rightarrow\infty$ follows similarly. Assume that $0\le q\le s-1$.  We compute the convolution explicitly using Stone's decomposition.
\begin{align}
\Phi*z(x)=&\int_0^x z(x-y)\Phi_2(\d y)+\int_0^x z(x-y)\varphi_1(y)\d y.
\label{}
\end{align}
The first integral can be split as
\begin{align}
\int_0^x z(x-y)\Phi_2(\d y)=\int_0^{x/2} z(x-y)\Phi_2(\d y)+\int_{x/2}^x z(x-y)\Phi_2(\d y),
\end{align}
where
\begin{align}
&\absol{\int_0^{x/2} z(x-y)\Phi_2(\d y)}\le\supre{}{y^\prime\ge x/2}\absol{z(y^\prime)}\norm{\Phi_2}_{\text{t.v.}}=o(x^{-r})
\label{}\\
&\absol{\int_{x/2}^x z(x-y)\Phi_2(\d y)}\le\norm{z}_{\infty}\Phi_2((x/2,\infty))=o(x^{-q}).
\end{align}
The second integral can be written as
\begin{align}
\int_0^x z(y)\varphi_1(x-y)\d y=&\int_{0}^{x/2}z(y)\varphi_1(x-y)\d y+\int_{x/2}^{x}z(y)\varphi_1(x-y)\d y.
\end{align}

From Lemma \ref{DecayRatesLemma}, we have $\tilde{\varphi}_1(x):=\varphi_1(x)-m=o(x^{-q})$ as $x\rightarrow\infty$. This entails,
\begin{align}
\absol{\int_{0}^{x/2}z(y)\tilde{\varphi}_1(x-y)\d y}\le&\int_{0}^{x/2}z(y)\absol{\tilde{\varphi}_1(x-y)}\d y
\label{}\\
\le& \supre{}{y^\prime\ge x/2}\absol{\tilde{\varphi}_1(y^\prime)}\int_{0}^\infty z(y)\d y
\label{}\\
=&o(x^{-q})
\label{}
\end{align}
On the other hand,
\begin{align}
\absol{\int_{x/2}^{x}z(y)\tilde{\varphi}_1(x-y)\d y}\le&\int_{x/2}^{x}z(y)\absol{\tilde{\varphi}_1(x-y)}\d y
\label{}\\
\le& \norm{\tilde{\varphi}_1}_{\infty}\int_{x/2}^{x}z(y)\d y,
\label{}
\end{align}
Let $\epsilon>0$ be arbitrary. Since $z(x)=o(x^{-r})$ as $x\rightarrow\infty$, there exists $x_0>0$ such that $x^r\absol{z(x)}<\epsilon$ for all $x>x_0$. Then, for all $x>2x_0$, 
\begin{align}
\int_{x/2}^{x}\absol{z(y)}\d y<&\epsilon\int_{x/2}^{x}y^{-r}\d y
\label{}\\
<&\epsilon\frac{(x/2)^{1-r}}{r-1},
\end{align}
from which we obtain that
\begin{align}
\int_{x/2}^{x}\absol{z(y)}\d y=o(x^{1-r}).
\end{align}
This means that,
\begin{align}
&\int_{0}^{x/2}z(y)\varphi_1(x-y)\d y=m\int_0^{x/2}z(y)\d y+o(x^{-q}),
\label{}\\
&\int_{x/2}^{x}z(y)\varphi_1(x-y)\d y=m\int_{x/2}^{x}z(y)\d y+o(x^{1-r}).
\end{align}
Ultimately, we get,
\begin{align}
\Phi*z(x)=&m\int_0^\infty z(y)\d y+o(x^{-q})+o(x^{-r})+o(x^{1-r})+o(x^{-q})
\label{}\\
=&m\int_0^\infty z(y)\d y+o(x^{\max\cbra{1-r,-q}}),
\end{align}
which is the desired result.
}

\section{Regenerative processes}
\label{SectionRegenerative}
We now turn to the notion of \emph{regenerative processes}. Intuitively, such processes evolve in cycles between the epochs of a renewal process. Formally, we have the following definition.
\begin{Def}[\textbf{Sec VI. p.169 in \cite{Asmussen2003}}]
\label{DefRegenerative}
Let $\cbra{X_t}_{t\ge0}$ be a stochastic process with state space $E$. We call $\cbra{X_t}_{t\ge0}$ \emph{regenerative} (pure or delayed) if there exists a renewal process (pure or delayed) $\cbra{\hat{S}_n}=\cbra{\tau_0+\tau_1+\dots+\tau_n}$ such that for each $n\ge0$, the \emph{post-$\hat{S}_n$ process}
\begin{align}
\theta_{\hat{S}_n}X:=\rbra{\tau_{n+1},\tau_{n+2},\dots,\cbra{X_{\hat{S}_n+t}}_{t\ge0}},
\end{align}
is independent of $\tau_0,\dots,\tau_n$, and its distribution does not depend upon $n$.
\end{Def}
Since the renewals are stopping times for a renewal process, from the strong Markov property it is obvious that the forward and backward recurrence times from Definition \ref{RecurrenceTimes} are regenerative. We can show as well that the compensator of a zero-delayed renewal process that satisfies assumption \textbf{(A1)} can be written in terms of a regenerative process. Let $\mu:[0,\infty)\rightarrow[0,\infty)$ be the nonnegative measurable function given by
\begin{align}
\mu(x)=\frac{f(x)}{1-\int_0^xf(y)\d y},\quad x\ge0.
\end{align}
From \cite[Sec. II Theorem T7]{BremaudQueues}, we have that the process $\cbra{M(t)}_{t\ge0}$ given as
\begin{align}
M(t)=N(t)-\int_0^t\mu(u-S_{N(u)})\d u,\quad t\ge0,
\end{align}
is an $(\cF_t^N)$-martingale, with \emph{compensator} $\Lambda$: 
\begin{align}
\Lambda(t):=\int_0^t\mu(u-S_{N(u)})\d u,\quad t\ge0.
\label{CompensatorDef}
\end{align}
Notice that $\Lambda$ can be rewritten as
\begin{align}
\Lambda(t)=&\sum_{i=1}^{N(t)-1}\xi_i+\int_{S_{N(t)-1}}^{t}\mu(u-S_{N(t)-1})\d u
\label{}\\
:=&\sum_{i=1}^{N(t)-1}\int_{S_{i-1}}^{S_{i}}\mu(u-S_{i-1})\d u+\int_{S_{N(t)-1}}^{t}\mu(u-S_{N(t)-1})\d u,
\end{align}
where the random variables $\xi_i$ are independent and identically distributed with
\begin{align}
\xi_i\dist\int_0^{\tau_1}\mu(u)\d u\quad i=1,2\dots
\end{align}
Observe that the process
\begin{align}
\int_{S_{N(t)-1}}^{t}\mu(u-S_{N(t)-1})\d u
\end{align}
is regenerative. In effect, let $n\ge0$ be arbitrary, then from the strong Markov property of the backward recurrence time $A_t$ we have
\begin{align}
\int_{S_{N(t+S_n)-1}}^{t+S_n}\mu(u-S_{N(t+S_n)-1})\d u=&\int_{0}^{t+S_n-S_{N(t+S_n)-1}}\mu(v)\d v
\label{}\\
=&\int_{0}^{A_{t+S_n}}\mu(v)\d v
\label{}\\
\dist& \int_{0}^{A_{t}}\mu(v)\d v
\end{align}
and since this is a measurable function whose argument is the regenerative process $A_t$, from \cite[Proposition VI.1.1]{Asmussen2003}, we can conclude that this process is itself regenerative.

Now, let $\cbra{X_t}_{t\ge0}$ be a regenerative process. It will be useful to compute the distribution of the maximum of $\bar{X}_T:=\max_{0\le t\le T}X_t$ for any $T>0$. To do this, we use the following result of Rootzén, which we state without proof. 
\begin{Thm}[\textbf{Rootz\'en \cite{RootzenRegenerative}}] 
\label{MaximumDistributionAsmussen}
Let $\cbra{X_t}_{t\ge0}$ be a regenerative process whose renewal process satisfies assumption \textbf{(B0)}. Assume that the state space of $X$ is a real interval and define 
\begin{align}
G(x):=\bP_0(\bar{X}_{\tau_1}\le x),\quad F_T(x):=\bP(\bar{X}_T\le x).
\end{align}
Then, when $G$ does not have finite support,
\begin{align}
\lim_{T\rightarrow\infty}\supre{}{x\ge0}\absol{F_T(x)-G(x)^{mT}}=0.
\end{align}
\end{Thm}

\section{Convergence of processes}
\label{SectionConvergence}
Now we can proceed with the proof of our Theorems \ref{MaxCompensator} and \ref{MaxRecurrenceTime}. We work in the context of regenerative processes introduced in the previous section.
\Proof[Proof of Theorem \ref{MaxCompensator}]{As we noted above, for all $i=1,2,\dots$, 
\begin{align}
\xi_i:=\int_{S_{i-1}}^{S_{i}}\mu(s-S_{{i-1}})\d s\dist\int_0^{\tau_1}\mu(s)\d s.
\end{align}
We now look for the distribution of the latter. Let $x>0$, we have
\begin{align}
\mu(x)=\frac{f(x)}{1-F(x)}=-\frac{\d}{\d x}\log(1-F(x)),
\end{align}
or in other words,
\begin{align}
\exp\rbra{-\int_0^x\mu(s)\d s}=1-F(x),
\end{align}
which implies that
\begin{align}
G(x):=&\bP\rbra{\int_0^{\tau_1}\mu(s)\d s\le x}
\label{}\\
=&\bP\rbra{\exp\rbra{-\int_0^{\tau_1}\mu(s)\d s}\ge e^{-x}}
\label{}\\
=&\bP\rbra{F(\tau_1)\le 1-e^{-x}}=1-e^{-x},
\end{align}
where in the last equality we used that if $F(x)=\bP(\tau_1\le x)$, then the random variable $F(\tau_1)\sim\text{Uniform}[0,1]$. This means that
\begin{align}
&\int_0^{\tau_1}\mu(s)\d s\sim\exp(1).
\label{}
\end{align}
We have for any $p>0$ that
\begin{align}
\supre{}{v\in[0,1]}\frac{1}{T^p}\sbra{\Lambda(Tv)-\Lambda(S_{N(Tv)-1})}=&\supre{}{v\in[0,1]}\frac{1}{T^p}\int_{S_{N(Tv)-1}}^{Tv}\mu(s-S_{N(Tv)-1})\d s
\label{}\\
\le &\frac{1}{T^p}\maxi{}{k\le N(T)}\xi_k,
\label{}
\end{align}
and since $G$ does not have finite support, from Theorem \ref{MaximumDistributionAsmussen}, we know that $F_{T}(x):=\bP\rbra{\maxi{}{k\le N(T)}\xi_k\le x}$ can be uniformly approximated by $G(x)^{mT}$, in the sense that
\begin{align}
\lim_{T\rightarrow\infty}\supre{}{x\ge0}\absol{F_{T}(x)-G(x)^{mT}}=0,
\end{align}
where $G(x)=\bP\rbra{\int_0^{\tau_1}\mu(s)\d s\le x}=1-e^{-x}$. If we take $x=\epsilon T^p$ for $\epsilon>0$, we have,
\begin{align}
G\rbra{\epsilon T^p}^{mT}=\rbra{1-e^{-\epsilon T^p}}^{mT}=&\exp\cbra{mT\log\rbra{1-e^{-\epsilon T^p}}}
\label{}\\
=&\exp\cbra{mT\rbra{e^{-\epsilon T^p}(1+o(1))}}\tend{}{T\rightarrow\infty}1,
\end{align}
which proves that for any $p>0$
\begin{align}
&\lim_{T\rightarrow\infty}\bP\rbra{\frac{1}{T^p}\maxi{}{k\le N(T)}\xi_k\le\epsilon}=1
\label{}
\end{align}
and hence,
\begin{align}
\rbra{\frac{1}{T^p}\int_{S_{{N(Tv)-1}}}^{Tv}\mu(s-S_{{N_R(Tv)-1}})\d s}_{v\in[0,1]}\tend{d}{T\rightarrow\infty}0.
\end{align}
The proof is complete.
}
We proceed similarly in the case of the recurrence times. As we will see, the dominating process will be the maximum process of the interarrival times.
\Proof[Proof of Theorem \ref{MaxRecurrenceTime}]{
Let $p\le s$. Notice that 
\begin{align}
&\supre{}{v\in[0,1]}\frac{1}{T^{1/p}}\rbra{Tv-S_{{N(Tv)-1}}}\le \frac{1}{T^{1/p}}\maxi{}{k\le N(T)}\tau_k,
\label{}\\
&\supre{}{v\in[0,1]}\frac{1}{T^{1/p}}\rbra{S_{{N(Tv)}}-Tv}\le \frac{1}{T^{1/p}}\maxi{}{k\le N(T)}\tau_k.
\label{}
\end{align}
It is clear that if $F(x):=\bP(\tau_1\le x)$ has finite support, the result is trivial. Thus, we assume that $F$ does not have a finite support. In this case we can approximate the distribution of the maximum by
\begin{align}
F(x)^{mT},
\end{align}
and by taking $x=\epsilon T^{1/p}$ we can compute,
\begin{align}
F(\epsilon T^{1/p})^{mT}=&\exp\cbra{mT\log\rbra{F(\epsilon T^{1/p})}}
\label{}\\
=&\exp\cbra{mT\log\rbra{1-\bar{F}(\epsilon T^{1/p})}}
\label{}\\
=&\exp\cbra{mT\rbra{\bar{F}(\epsilon T^{1/p})(1+o(1))}}.
\end{align}
Since $\bE\sbra{\tau_1^s}<\infty$, it follows by the DCT that
\begin{align}
&x^p\bar{F}(x)=x^p\bP\rbra{\tau_1>x}\le\bE\sbra{\tau_1^p\; ;\; \tau_1>x}\tend{}{x\rightarrow\infty}0.
\label{}
\end{align}
The previous reasoning implies that
\begin{align}
&T\bar{F}(\epsilon T^{1/p})\tend{}{T\rightarrow\infty}0,
\label{}\\
&\exp\cbra{mT\rbra{\bar{F}(\epsilon T^{1/p})(1+o(1))}}\tend{}{T\rightarrow\infty}1,
\label{}\\
&\lim_{T\rightarrow\infty}\bP\rbra{\frac{1}{T^{1/p}}\maxi{}{k\le N(T)}\tau_k\le\epsilon}=1,
\end{align}
from which we can conclude that
\begin{align}
\rbra{\frac{1}{T^{1/p}}\rbra{Tv-S_{{N(Tv)-1}}},\frac{1}{T^{1/p}}\rbra{S_{{N(Tv)}}-Tv}}_{v\in[0,1]}\tend{d}{T\rightarrow\infty}(0,0),
\end{align}
which completes the proof.
}
\section{Appendix: Some results from the general theory of renewal processes}
\label{SecctionAppendixA}
The proofs of the results stated in this section are included for the completeness of the paper.
\begin{Lem}
\label{RHPALinearSolution}
Assume that \textbf{(B0)} holds. Then, the renewal equation
\begin{align}
\label{AppendixAEquation1}
Z(t)=z(t)+\int_0^tZ(t-x)F(\d x),
\end{align}
has the linear solution $Z_0(t)=mt$, $t\ge0$, if and only if
\begin{align}
z(t)=m\int_0^t\bar{F}(x)\d x,
\end{align}
where $\bar{F}=1-F$.
\end{Lem}
\Proof{Let us suppose that equation \eqref{AppendixAEquation1} has a linear solution $Z_0(t)=mt$, which entails
\begin{align}
mt=z(t)+\int_0^tm(t-x)F(\d x).
\end{align}
Solving for $z(t)$ we get
\begin{align}
z(t)=&mt-\int_0^t m(t-x)F(\d x)
\label{}\\
=&mt-\int_0^tm\int_x^t\d u F(\d x)
\label{}\\
=&mt-\int_0^tm\int_0^uF(\d x)\d u 
\label{}\\
=&\int_0^tm\d x-\int_0^tmF(x)\d x
\label{}\\
=&\int_0^tm(1-F(x))\d x
\label{}\\
=&\int_0^tm\bar{F}(x)\d x.
\end{align}
Reversing the order of the above steps we can see that the function $z(t)=\int_0^tm\bar{F}(x)\d x$ generates a renewal Equation with solution
\begin{align}
mt=Z_0(t)=\int_0^tm\int_0^{t-u}\bar{F}(x)\d x\Phi(\d u).
\end{align}
The proof is now complete.
}
The distribution of the forward recurrence time can be found through a renewal argument as follows
\begin{Lem}
\label{DistBt}
The distribution function of the forward recurrence time is given for $x\ge0$ as
\begin{align}
\bP\rbra{B_t\le x}=\int_0^t F\rbra{(t-u,t+x-u]}\Phi(\d u).
\end{align}
\end{Lem}
\Proof{Write $Z_x(t)=\bP(B_t\le x)$, and $z_x(t)=\bP(B_t\le x; \tau_1>t)$. Then,
\begin{align}
Z_x(t)=&z_x(t)+\bP(B_t\le x; \tau_1\le t)
\label{}\\
=&z_x(t)+\int_0^t\bP\rbra{B_t\le x\mid \tau_1=u}F(\d u)
\label{}\\
=&z_x(t)+\int_0^t\bP\rbra{B_{t-u}\le x}F(\d u)\quad\text{(by the Strong Markov property of $B_t$)}
\label{}\\
=&z_x(t)+\int_0^tZ_x(t-u)F(\d u)=\int_0^t z_x(t-u)\Phi(\d u).
\end{align}
Now, analyzing $z_x(t)$, we have,
\begin{align}
z_x(t)=&\bP(B_t\le x; \tau_1>t)
\label{}\\
=&\bP(\tau_1-t\le x; \tau_1>t)
\label{}\\
=&F(t+x)-F(t).
\end{align}
Therefore, we obtain,
\begin{align}
Z_x(t)=\int_0^t F\rbra{(t-u,t+x-u]}\Phi(\d u),
\end{align}
which concludes the proof.
}
\section{Acknowledgements}
To finalize this paper, I want to sincerely thank my PhD supervisor, Dr. Kouji Yano, for his continuous support throughout the writing of this work. The several readings he gave to the manuscript and all his insightful comments made possible the coming together of this paper.

\bibliographystyle{plain}

\begin{thebibliography}{10}
\bibitem{Asmussen2003}
S.~Asmussen.
\newblock {\em Applied Probability and Queues}.
\newblock Applications of mathematics : stochastic modelling and applied
  probability. Springer, 2003.

\bibitem{asmussen2013asymptotics}
S.~Asmussen, S.~Foss, and D.~Korshunov.
\newblock Asymptotics for sums of random variables with local subexponential
  behaviour, 2013.

\bibitem{AthreyaNey}
K.~B. Athreya and P.~Ney.
\newblock A new approach to the limit theory of recurrent markov chains.
\newblock {\em Transactions of the American Mathematical Society},
  245:493--501, 1978.

\bibitem{DVJ}
D.~J. Daley and D.~Vere-Jones.
\newblock {\em An introduction to the theory of point processes. {V}ol. {I}}.
\newblock Probability and its Applications (New York). Springer-Verlag, New
  York, second edition, 2003.
\newblock Elementary theory and methods.

\bibitem{GakisRecurrence}
K.G. Gakis and B.D. Sivazlian.
\newblock The use of multiple integrals in the study of the backward and
  forward recurrence times for the ordinary renewal process.
\newblock {\em Stochastic Analysis and Applications}, 10(4):409--416, 1992.

\bibitem{godreche2000statistics}
C.~Godreche and J.~M. Luck.
\newblock Statistics of the occupation time of renewal processes.
\newblock {\em Journal of Statistical Physics}, 104:489--524, 2001.

\bibitem{LindvallCoupling}
Torgny Lindvall.
\newblock On coupling of continuous-time renewal processes.
\newblock {\em Journal of Applied Probability}, 19(1):82--89, 1982.

\bibitem{LundExponential}
Robert~B. Lund, Sean~P. Meyn, and Richard~L. Tweedie.
\newblock Computable exponential convergence rates for stochastically ordered
  markov processes.
\newblock {\em The Annals of Applied Probability}, 6(1):218--237, 1996.

\bibitem{BremaudQueues}
Brémaud Pierre.
\newblock {\em Point processes and queues : martingale dynamics / Pierre
  Brémaud}.
\newblock Springer series in statistics. Springer-Verlag, New York, cop. 1981.

\bibitem{RootzenRegenerative}
Holger Rootzén.
\newblock Maxima and exceedances of stationary markov chains.
\newblock {\em Advances in Applied Probability}, 20(2):371--390, 1988.

\bibitem{Sgibnev1981}
Mikhail~Sergeyevich Sgibnev.
\newblock On the renewal theorem in the case of infinite variance.
\newblock {\em Sibirskii Matematicheskii Zhurnal}, 22(5):178--189, 1981.

\bibitem{Stone1966}
Charles Stone.
\newblock {On Absolutely Continuous Components and Renewal Theory}.
\newblock {\em The Annals of Mathematical Statistics}, 37(1):271 -- 275, 1966.

\bibitem{thorisson2000coupling}
H.~Thorisson.
\newblock {\em Coupling, Stationarity, and Regeneration}.
\newblock Probability and Its Applications. Springer New York, 2000.

\bibitem{WillmotLundberg}
Gordon~E. Willmot, Jun Cai, and X.~Sheldon Lin.
\newblock Lundberg inequalities for renewal equations.
\newblock {\em Advances in Applied Probability}, 33(3):674--689, 2001.

\bibitem{ChuancunNonExpAsym}
Chuancun Yin and Junsheng Zhao.
\newblock Nonexponential asymptotics for the solutions of renewal equations,
  with applications.
\newblock {\em Journal of Applied Probability}, 43(3):815--824, 2006.

\end{thebibliography}

\end{document}